\documentclass[leqno,12pt]{amsart}
\usepackage{amsfonts}
\usepackage{amsmath,amssymb,amsthm}
\newcommand{\R}{\mathbb R}

\newcommand{\E}{\mathbb E}

\renewcommand{\span}{\mathrm{span}}
\newcommand{\tr}{\mathrm{tr}}

\newtheorem{thm}{Theorem}[section]

\theoremstyle{definition}

\theoremstyle{remark}

\newcommand{\ds}{\displaystyle}

\begin{document}

\title[Constant Mean Curvature Rotational Surfaces]
{Rotational Surfaces with Constant Mean Curvature in Pseudo-Euclidean
4-Space with Neutral Metric}

\author{Yana  Aleksieva,  Velichka Milousheva}

\address{Faculty of Mathematics and Informatics, Sofia University,
5 James Bourchier blvd., 1164 Sofia, Bulgaria}
\email{yana\_a\_n@fmi.uni-sofia.bg}
\address{Institute of Mathematics and Informatics, Bulgarian Academy of Sciences,
Acad. G. Bonchev Str. bl. 8, 1113, Sofia, Bulgaria;   "L.
Karavelov" Civil Engineering Higher School, 175 Suhodolska Str., 1373 Sofia, Bulgaria}
\email{vmil@math.bas.bg}


\begin{abstract}
We give the classification of constant mean curvature  rotational  
surfaces of elliptic, hyperbolic, and parabolic type in the four-dimensional pseudo-Euclidean space with neutral metric.
\end{abstract}

\maketitle

\section{Introduction}

In the Minkowski 4-space $\E^4_1$ there exist three
types of rotational surfaces with two-dimensional axis -- rotational
surfaces of elliptic, hyperbolic or parabolic type, known also as
surfaces invariant under spacelike rotations, hyperbolic rotations
or screw rotations,  respectively.
A \emph{rotational surface of elliptic type} is an orbit of a regular
curve under the action of the orthogonal transformations of $\E^4_1$
which leave a timelike plane point-wise fixed. Similarly, a \emph{rotational
surface of hyperbolic type} is an orbit of a  regular curve under
the action of the orthogonal transformations of $\E^4_1$ which leave
a spacelike plane point-wise fixed. A \emph{rotational surface of
parabolic type} is an an orbit of a  regular curve under the action
of the orthogonal transformations of $\E^4_1$ which leave a
degenerate plane point-wise fixed.

The marginally trapped surfaces in Minkowski 4-space which are
invariant under spacelike rotations (rotational surfaces of elliptic
type) were classified by S. Haesen and  M. Ortega  in
\cite{Haesen-Ort-2}.   The classification of marginally trapped
surfaces   in $\E^4_1$ which are invariant under  boost
transformations (rotational surfaces of hyperbolic type) was
obtained  in \cite{Haesen-Ort-1} and the classification of
marginally trapped surfaces which are invariant under  screw
rotations (rotational surfaces of parabolic type) is given in
\cite{Haesen-Ort-3}.

Motivated by the classification results of  S. Haesen and M. Ortega about
marginally trapped rotational surfaces in the Minkowski space,
in \cite{GM5} G. Ganchev and the second author considered three types of rotational surfaces
in the four-dimensional  pseudo-Euclidean space $\E^4_2$, namely
rotational surfaces of elliptic, hyperbolic, and parabolic type,
which are analogous to the three types of rotational surfaces in  $\E^4_1$. They classified  all quasi-minimal rotational surfaces of elliptic, hyperbolic, and  parabolic type.

Constant mean curvature surfaces in arbitrary spacetime
are important objects for the special role they play in the theory of
general relativity. The study of constant mean curvature surfaces
(CMC surfaces) involves not only geometric methods but also PDE
and complex analysis, that is why the theory of CMC surfaces is of
great interest not only for mathematicians but also for physicists
and engineers. Surfaces with constant mean curvature in
Minkowski space have been studied intensively in  the last years. See for example \cite{Bran}, \cite{Chav-Can}, \cite{Liu-Liu-1}, \cite{Lop-2}, \cite{Sa}.

Classification results for rotational surfaces in three-dimensional
space forms satisfying some classical extra conditions have also been obtained.
For example, a classification of all timelike and
spacelike hyperbolic rotational surfaces with non-zero constant mean
curvature in the three-dimensional de Sitter space $\mathbb{S}^3_1$
is given in \cite{Liu-Liu} and a classification of the spacelike and
timelike Weingarten rotational surfaces of the three types in
$\mathbb{S}^3_1$ is found in \cite{Liu-Liu-2}. Chen spacelike rotational surfaces of hyperbolic or
elliptic type are described in  \cite{GM3}.

In the present paper we study Lorentz  rotational surfaces of elliptic,  hyperbolic, and parabolic type and give the classification of all such surfaces  with non-zero  constant mean curvature in $\E^4_2$.

\section{Preliminaries}

Let  $\E^4_2$ be the pseudo-Euclidean 4-space endowed with the canonical pseudo-Euclidean metric of index 2 given by
$g_0 = dx_1^2 + dx_2^2 - dx_3^2 - dx_4^2,$
where $(x_1, x_2, x_3, x_4)$ is a rectangular coordinate system of $\E^4_2$. As usual, we denote by
$\langle \, , \rangle$ the indefinite inner scalar product with respect to $g_0$.
 A non-zero vector $v$ is called  \emph{spacelike} (respectively, \emph{timelike}) if $\langle v, v \rangle > 0$ (respectively, $\langle v, v \rangle < 0$).
 A vector $v$ is called \emph{lightlike} if it is nonzero and satisfies $\langle v, v \rangle = 0$.

A surface $M^2_1$ in $\E^4_2$ is called \emph{Lorentz}  if the
induced  metric $g$ on $M^2_1$ is Lorentzian, i.e. at each point $p
\in M^2_1$ we have the following decomposition
$\E^4_2 = T_pM^2_1 \oplus N_pM^2_1$
with the property that the restriction of the metric
onto the tangent space $T_pM^2_1$ is of
signature $(1,1)$, and the restriction of the metric onto the normal space $N_pM^2_1$ is of signature $(1,1)$.

Denote by $\nabla$ and $\nabla'$ the Levi Civita connections of $M^2_1$  and $\E^4_2$, respectively.
Let $x$ and $y$ denote vector fields tangent to $M^2_1$ and $\xi$ be a normal vector field.
The formulas of Gauss and Weingarten are given respectively by
$$\begin{array}{l}
\vspace{2mm}
\nabla'_xy = \nabla_xy + \sigma(x,y);\\
\vspace{2mm}
\nabla'_x \xi = - A_{\xi} x + D_x \xi,
\end{array}$$
where $\sigma$ is the second fundamental form, $D$ is the normal
connection, and $A_{\xi}$ is the shape operator  with respect to
$\xi$. In general, $A_{\xi}$ is not diagonalizable.

The mean curvature vector  field $H$ of the surface $M^2_1$
is defined as $H = \frac{1}{2}\,  \tr\, \sigma$.
A  surface $M^2_1$  is called \emph{minimal} if its mean curvature vector vanishes identically, i.e. 
$H =0$.
A  surface $M^2_1$  is \emph{quasi-minimal} if its
mean curvature vector is lightlike at each point, i.e. $H \neq 0$ and $\langle H, H \rangle =0$.
In this paper we consider Lorentz surfaces in $\E^4_2$ for which  $\langle H, H \rangle = const \neq 0$.

\section{Lorentz  rotational surfaces with constant mean curvature in $\E^4_2$}

 Let $Oe_1e_2e_3e_4$ be a fixed orthonormal coordinate system in the pseudo-Euclidean space $\E^4_2$ such that $\langle e_1, e_1 \rangle =
\langle e_2, e_2 \rangle  = 1, \, \langle e_3, e_3 \rangle =  \langle e_4, e_4 \rangle = -1$. Lorentz rotational surfaces of elliptic, hyperbolic, and parabolic type are defined in \cite{GM5}. Here we shall present shortly the construction.

First we consider rotational surfaces of elliptic type.
Let $c: \widetilde{z} = \widetilde{z}(u), \,\, u \in J$ be a smooth spacelike curve lying in the three-dimensional subspace $\E^3_1 = \span\{e_1, e_2, e_3\}$ and
parameterized by
$\widetilde{z}(u) = \left( x_1(u), x_2(u),  r(u), 0 \right); \; u \in J$.
Without loss of generality we assume that $c$ is
parameterized by the arc-length, i.e. $(x_1')^2 + (x_2')^2 - (r')^2 = 1$, and  $r(u)>0, \,\, u \in J$.

Let  $\mathcal{M}'$ be the surface in $\E^4_2$ defined by
\begin{equation} \label{E:Eq-1}
\mathcal{M}': z(u,v) = \left( x_1(u), x_2(u), r(u) \cos v, r(u) \sin v\right);
\quad u \in J,\,\,  v \in [0; 2\pi).
\end{equation}
The surface $\mathcal{M}'$, defined by \eqref{E:Eq-1}, is a Lorentz surface  in $\E^4_2$, obtained by the rotation of the spacelike curve $c$ about
the two-dimensional Euclidean plane $Oe_1e_2$. It is called  a \emph{rotational surface of elliptic type}.

One can also obtain a rotational surface of elliptic type in  $\E^4_2$ using rotation of  a timelike curve about the two-dimensional plane $Oe_3e_4$ (see \cite{GM5}).

\vskip 2mm
Next, we consider rotational surfaces of hyperbolic type.
Let $c: \widetilde{z} = \widetilde{z}(u), \,\, u \in J$ be a smooth spacelike curve, lying in the three-dimensional subspace $\E^3_1 = \span\{e_1, e_2, e_4\}$ of $\E^4_2$
and parameterized by
$\widetilde{z}(u) = \left( r(u), x_2(u), 0, x_4(u) \right); \; u \in J$.
Without loss of generality we assume that $c$ is
parameterized by the arc-length, i.e. $(r')^2 +(x_2')^2 - (x_4')^2  = 1$, and   $r(u)>0, \,\, u \in J$.

Let $\mathcal{M}''$ be the surface in $\E^4_2$ defined by
\begin{equation} \label{E:Eq-2}
\mathcal{M}'': z(u,v) = \left(r(u) \cosh v,  x_2(u), r(u) \sinh v, x_4(u) \right);
\quad u \in J,\,\,  v \in \R.
\end{equation}
The surface $\mathcal{M}''$, defined by \eqref{E:Eq-2}, is a Lorentz surface  in $\E^4_2$, obtained by hyperbolic rotation of the spacelike curve $c$ about
the two-dimensional Lorentz plane $Oe_2e_4$.
$\mathcal{M}''$ is called a \emph{rotational surface of hyperbolic type}.

Similarly, one can obtain a rotational surface of hyperbolic type  using hyperbolic rotation of  a timelike curve lying in $\span\{e_2, e_3, e_4\}$ about the two-dimensional Lorentz plane $Oe_2e_4$.
Rotational surfaces of hyperbolic type can also be obtained by  hyperbolic rotations of spacelike or timelike curves about the two-dimensional Lorentz planes $Oe_1e_3$,
 $Oe_1e_4$ and $Oe_2e_3$. 

\vskip 2mm
Now, we shall consider rotational surfaces of parabolic type in $\E^4_2$.
For convenience, in the parabolic case we  use the pseudo-orthonormal base $\{e_1,
e_4, \xi_1, \xi_2 \}$  of $\E^4_2$, such that $ \xi_1=\ds{
\frac{e_2 + e_3}{\sqrt{2}},\,\, \xi_2= \frac{ - e_2 +
e_3}{\sqrt{2}}}$. Note that
$\langle \xi_1, \xi_1 \rangle =0; \; \langle \xi_2, \xi_2 \rangle =0; \; \langle \xi_1, \xi_2 \rangle = -1.$

Let $c$ be a spacelike curve lying in the subspace $\E^3_1 = \span\{e_1, e_2, e_3\}$ of $\E^4_2$ and parameterized by
$\widetilde{z}(u) =  x_1(u)\, e_1 + x_2(u) \, e_2 + x_3(u) \, e_3; \; u \in J$,
or equivalently,
$$ \widetilde{z}(u) =  x_1(u)\, e_1 + \frac{x_2(u) + x_3(u)}{\sqrt{2}} \, \xi_1 + \frac{ - x_2(u) + x_3(u)}{\sqrt{2}} \, \xi_2; \quad u \in J.$$
Denote $f(u) =  \ds{\frac{x_2(u) + x_3(u)}{\sqrt{2}}}$, $g(u) = \ds{\frac{ - x_2(u) + x_3(u)}{\sqrt{2}}}$. Then
$ \widetilde{z}(u) =  x_1(u)\, e_1 + f(u) \, \xi_1 + g(u) \, \xi_2.$
Without loss of generality we assume that $c$ is
parameterized by the arc-length, i.e. $(x_1')^2 +(x_2')^2 - (x_3')^2  = 1$, or equivalently  $(x_1')^2 - 2f'g' = 1$.

A rotational surface of parabolic type is defined 
in the following way:
\begin{equation} \label{E:Eq-3}
\mathcal{M}''': z(u,v) = x_1(u)\, e_1 + f(u) \, \xi_1 + (-v^2 f(u) + g(u)) \, \xi_2 + \sqrt{2}\, v f(u) \, e_4;
\; u \in J,\,\,  v \in \R.
\end{equation}
The rotational axis is the  two-dimensional plane spanned by $e_1$ (a spacelike vector field) and $\xi_1$ (a lightlike vector field).

Similarly, one can obtain a rotational surface of parabolic type using a timelike curve lying in the subspace $\span\{e_2, e_3, e_4\}$ (see \cite{GM5}).

In what follows, we find all CMC Lorentz rotational surfaces of elliptic, hyperbolic, and parabolic type.

\vskip 3mm
\subsection{Constant mean curvature rotational surfaces of elliptic type}

Let us  consider the surface $\mathcal{M}'$ in $\E^4_2$, defined by \eqref{E:Eq-1}. 
The tangent space of $\mathcal{M}'$ is spanned by the vector fields
$z_u = \left(x_1', x_2', r' \cos v, r' \sin v  \right); \; 
 z_v = \left( 0, 0 , - r \sin v, r \cos v\right)$.
Hence, the coefficients of the first fundamental form of $\mathcal{M}'$ are
$E = \langle z_u, z_u \rangle = 1; \, F = \langle z_u, z_v \rangle = 0; \, G = \langle z_v, z_v \rangle = - r^2(u)$. Since the generating curve $c$ is a spacelike curve
parameterized by the arc-length, i.e. $(x_1')^2 + (x_2')^2 - (r')^2 = 1$, then $(x_1')^2 + (x_2')^2 = 1 + (r')^2 $ and $x_1' x_1'' + x_2 ' x_2'' = r' r''$.
We consider the  orthonormal tangent frame field
$X = z_u; \; Y = \ds{\frac{z_v}{r}},$
and the normal frame field $\{n_1, n_2\}$, defined by
\begin{equation} \label{E:Eq-5}
\begin{array}{l}
n_1 = \ds{\frac{1}{\sqrt{1+(r')^2}}\left(- x_2', x_1',0,0 \right)};\\
n_2 = \ds{\frac{1}{\sqrt{1+(r')^2}} \left(r' x_1', r' x_2', (1+(r')^2) \cos v, (1+(r')^2) \sin v \right)}.
\end{array}
\end{equation}
Note that $\langle X, X \rangle = 1; \; \langle X, Y \rangle = 0; \; \langle Y, Y \rangle = -1;$
$\langle n_1, n_1 \rangle = 1; \; \langle n_1, n_2 \rangle = 0; \; \langle n _2, n_2 \rangle = -1.$

Calculating the second partial derivatives of $z(u,v)$ and  the components of the second fundamental form, we obtain the formulas:
\begin{equation*} \label{E:Eq-7}
\begin{array}{l}
\vspace{1mm}
\sigma(X,X)= \ds{\frac{x_1' x_2'' - x_1'' x_2'}{\sqrt{1+(r')^2}} \, n_1 + \frac{r''}{\sqrt{1+(r')^2}}\, n_2},\\
\vspace{1mm}
\sigma(X,Y)= 0,\\
\vspace{1mm}
\sigma(Y,Y) =\qquad \qquad \qquad \quad \; - \ds{\frac{\sqrt{1+(r')^2}}{r}\, n_2},
\end{array}
\end{equation*}
which  imply that the normal mean curvature vector field $H$ of $\mathcal{M}'$ is expressed as follows
\begin{equation}  \label{E:Eq-8}
H = \frac{1}{2 r\sqrt{1+(r')^2}} \left( r(x_1' x_2'' - x_1'' x_2') \, n_1 + (r r'' + (r')^2 + 1) \, n_2 \right).
\end{equation}

Hence, $\langle H, H \rangle = \ds{\frac{r^2(x_1' x_2'' - x_1'' x_2')^2 - (r r'' + (r')^2 + 1) ^2} {4 r^2(1+(r')^2)}}$.
In the present paper we are interested in  rotational surfaces with non-zero constant mean curvature. So, we assume that $r^2(x_1' x_2'' - x_1'' x_2')^2 - (r r'' + (r')^2 + 1)^2 \neq 0$.

It follows from    \eqref{E:Eq-5} that
\begin{equation} \label{E:Eq-9}
\begin{array}{l}
\vspace{2mm}
\nabla'_X n_1 =- \ds{\frac{x_1' x_2'' - x_1'' x_2'}{\sqrt{1+(r')^2}} \, X + \frac{r'}{1+(r')^2}(x_1' x_2'' - x_1'' x_2')}\, n_2, \\
\vspace{2mm}
\nabla'_Y n_1 = 0,\\
\vspace{2mm}
\nabla'_X n_2 =\ds{\frac{r''}{\sqrt{1+(r')^2}} \, X + \frac{r'}{1+(r')^2}(x_1' x_2'' - x_1'' x_2') \, n_1}, \\
\vspace{2mm}
\nabla'_Y n_2 = \ds{\frac{\sqrt{1+(r')^2}}{r}\, Y}. 
\end{array}
\end{equation}

If  $x_1' x_2'' - x_1'' x_2' = 0$,\,  $r r'' + (r')^2 + 1 \neq 0$, then from  \eqref{E:Eq-9} we get
$\nabla'_X n_1 = \nabla'_Y n_1 = 0$,
 which imply that the normal vector field $n_1$ is constant. Hence,
the surface $\mathcal{M}'$ lies in the hyperplane $\E^3_2 =  \span \{X,Y,n_2\}$.

So, further we  consider rotational surfaces of elliptic  type satisfying $x_1' x_2'' - x_1'' x_2' \neq 0$ in an open interval $I \subset J$. 

In the next theorem we give a local description of all constant mean curvature rotational surfaces of elliptic  type.

\begin{thm}\label{T:cmc elliptic}
Given a smooth positive function $r(u): I \subset \R \rightarrow \R$, define the functions
$$\varphi(u) = \eta \int \ds{\frac{\sqrt{(r r'' + (r')^2 + 1)^2 \pm 4 C^2r^2(1+(r')^2)}}{r (1+(r')^2)}} \, du, \quad \eta = \pm 1, \; C = const \neq 0$$
and
\begin{equation} \notag
\begin{array}{l}
\vspace{2mm}
x_1(u) = \int  \sqrt{1+(r')^2} \,\cos \varphi(u) \, du,\\
\vspace{2mm}
x_2(u) =\int  \sqrt{1+(r')^2} \,\sin \varphi(u) \, du.
\end{array}
\end{equation}
Then the spacelike curve $c: \widetilde{z}(u) = \left( x_1(u), x_2(u),  r(u), 0 \right)$ is a generating curve of a
constant mean curvature rotational surface of elliptic  type.

Conversely, any constant mean curvature rotational surface of elliptic  type is locally constructed as above.
\end{thm}

\noindent
\emph{Proof:}
Let $\mathcal{M}'$ be a general  rotational surface of elliptic  type generated by a spacelike curve
$c: \widetilde{z}(u) = \left( x_1(u), x_2(u),  r(u), 0 \right);\,\, u \in J$.
We assume that $c$ is parameterized by the arc-length  and $x_1' x_2'' - x_1'' x_2' \neq 0$ for $u \in I \subset J$.
Using \eqref{E:Eq-8}  we get that $\mathcal{M}'$ is of constant mean curvature if and only if
 \begin{equation} \label{E:Eq-11}
\ds{\frac{r^2 (x_1' x_2'' - x_1'' x_2')^2 - (r r'' + (r')^2 + 1)^2}{4r^2(1+(r')^2)}}= \varepsilon C^2,\; \varepsilon = sign \langle H, H \rangle, \; C = const \neq 0.
\end{equation}
Since the curve $c$  is parameterized by the arc-length, we have   $(x_1')^2 + (x_2')^2 = 1 + (r')^2$, which implies that
there exists a smooth function $\varphi = \varphi(u)$ such that
\begin{equation} \label{E:Eq-12}
\begin{array}{l}
\vspace{2mm}
x_1'(u) = \sqrt{1+(r')^2} \,\cos \varphi(u),\\
\vspace{2mm}
x_2'(u) = \sqrt{1+(r')^2} \,\sin \varphi(u).
\end{array}
\end{equation}
It follows from  \eqref{E:Eq-12} that $x_1' x_2'' - x_1'' x_2' = (1+(r')^2) \varphi'$. Hence, condition \eqref{E:Eq-11} 
is written in terms of $r(u)$ and $\varphi(u)$ as follows:
\begin{equation} \label{E:Eq-13}
\varphi'(u) = \eta \ds{\frac{\sqrt{(r r'' + (r')^2 + 1)^2  \pm 4 C^2r^2(1+(r')^2)}}{r (1+(r')^2)}}, \quad \eta = \pm 1.
\end{equation}
Formula \eqref{E:Eq-13} allows us to recover the function $\varphi(u)$ from $r(u)$, up to integration constant.
Using formulas \eqref{E:Eq-12}, we can recover $x_1(u)$ and $x_2(u)$ from the functions $\varphi(u)$ and $r(u)$, up to integration constants.
Consequently, the constant mean curvature rotational surface of elliptic type $\mathcal{M}'$ is constructed as described in the theorem.

Conversely, given a smooth function $r(u) > 0$, we can define the function
$$\varphi(u) = \eta \int \ds{\frac{\sqrt{(r r'' + (r')^2 + 1)^2 \pm 4 C^2r^2(1+(r')^2)}}{r (1+(r')^2)}} \, du,$$
where $\eta = \pm 1$,  and consider the functions
\begin{equation} \notag
\begin{array}{l}
\vspace{2mm}
x_1(u) = \int  \sqrt{1+(r')^2} \,\cos \varphi(u) \, du,\\
\vspace{2mm}
x_2(u) =\int  \sqrt{1+(r')^2} \,\sin \varphi(u) \, du.
\end{array}
\end{equation}
A straightforward computation shows that the curve
$c: \widetilde{z}(u) = \left( x_1(u), x_2(u),  r(u), 0 \right)$ is a spacelike curve generating a constant mean curvature rotational surface of elliptic type according to formula \eqref{E:Eq-1}.

\qed

\vskip 2mm
\noindent
\textit{Remark}: 
In the special case when $r r'' + (r')^2 + 1 = 0$, i.e $r(u) = \pm \sqrt{-u^2+2au+b}, \; a=const\neq 0,\,  b=const$, 
the function $\varphi(u)$ is expressed as
$$\varphi(u)= \ds{\frac{2C}{\sqrt{a^2+b}}} \left(\frac{u-a}{2}\sqrt{-u^2+2au+b} + \frac{a^2+b}{2} \arcsin\frac{u-a}{\sqrt{a^2+b}} + d\right), \; d = const.$$

\vskip 3mm
\subsection{Constant mean curvature rotational surfaces of hyperbolic type}

Now, we shall consider the rotational surface of hyperbolic type  $\mathcal{M}''$, defined by \eqref{E:Eq-2}.
The tangent space of $\mathcal{M}''$ is spanned by the vector fields
$z_u = \left( r' \cosh v, x_2',  r' \sinh v, x_4' \right)$;
$z_v = \left(r \sinh v, 0, r \cosh v, 0 \right)$,
and the coefficients of the first fundamental form of $\mathcal{M}''$ are
$E = 1; \; F  = 0; \; G = - r^2(u)$.

The generating curve $c$ is a spacelike curve
parameterized by the arc-length, i.e. $(r')^2 + (x_2')^2 - (x_4')^2  = 1$, and hence $(x_4')^2 - (x_2')^2 = (r')^2 - 1 $.
We assume that $(r')^2 \neq 1$, otherwise the surface lies in a 2-dimensional plane.
Denote by $\varepsilon$ the sign of $(r')^2 - 1$.

We use the  orthonormal tangent frame field 
$X = z_u; \; Y = \ds{\frac{z_v}{r}}$
and the normal frame field $\{n_1, n_2\}$, defined by
\begin{equation} \label{E:Eq-5a}
\begin{array}{l}
\vspace{2mm}
n_1 = \ds{\frac{1}{\sqrt{\varepsilon((r')^2 - 1)}}\left(0, x_4', 0, x_2' \right)};\\
\vspace{2mm}
n_2 = \ds{\frac{1}{\sqrt{\varepsilon((r')^2 - 1)}} \left((1-(r')^2) \cosh v, - r' x_2',  (1-(r')^2) \sinh v, - r' x_4' \right)}.
\end{array}
\end{equation}
The frame field $\{X, Y, n_1, n_2\}$ satisfies
$\langle X, X \rangle = 1; \, \langle X, Y \rangle = 0; \, \langle Y, Y \rangle = -1; \, \langle n_1, n_1 \rangle = \varepsilon; \, \langle n_1, n_2 \rangle = 0; \, \langle n _2, n_2 \rangle = - \varepsilon$.

The mean curvature vector field $H$ of $\mathcal{M}''$ is expressed by the following formula
\begin{equation*}  \label{E:Eq-8a}
H = \ds{\frac{\varepsilon}{2 r\sqrt{\varepsilon((r')^2 - 1)}} \left( r(x_4' x_2'' - x_4'' x_2') \, n_1 - (r r'' + (r')^2 - 1) \, n_2 \right)}.
\end{equation*}
If $x_2' x_4'' - x_2'' x_4' = 0$,\,  $r r'' + (r')^2 -1 \neq 0$, $\mathcal{M}''$ lies in the hyperplane $\span \{X,Y,n_2\}$. So, we assume  that $x_2' x_4'' - x_2'' x_4' \neq 0$ in an open interval $I \subset J$. 

The local classification of constant mean curvature rotational surfaces of hyperbolic type is given by the following theorem:

\begin{thm}\label{T:cmc hyperbolic}
\emph{Case (A)}. Given a smooth positive function $r(u): I \subset \R \rightarrow \R$, such that $(r')^2 > 1$,  define the functions
$$\varphi(u) = \eta \int \ds{\frac{\sqrt{ (r r'' + (r')^2 - 1)^2 \pm 4 C^2r^2((r')^2-1)}}{r ((r')^2-1)}} \, du, \quad \eta = \pm 1, \; C = const \neq 0,$$
and
\begin{equation} \notag
\begin{array}{l}
\vspace{2mm}
x_2(u) = \int  \sqrt{(r')^2-1} \,\sinh \varphi(u) \, du,\\
\vspace{2mm}
x_4(u) =\int  \sqrt{(r')^2-1} \,\cosh \varphi(u) \, du.
\end{array}
\end{equation}
Then the spacelike curve $c: \widetilde{z}(u) = \left( r(u), x_2(u), 0, x_4(u) \right)$ is a generating curve of a
constant mean curvature rotational surface of hyperbolic type. 

\emph{Case (B)}. Given a smooth positive function $r(u): I \subset \R \rightarrow \R$,  such that $(r')^2 < 1$, define the functions
$$\varphi(u) = \eta \int \ds{\frac{\sqrt{ (r r'' + (r')^2 - 1)^2 \pm 4 C^2r^2((r')^2-1)}}{r ((r')^2-1)}} \, du, \quad \eta = \pm 1, \; C = const \neq 0,$$
and
\begin{equation} \notag
\begin{array}{l}
\vspace{2mm}
x_2(u) = \int  \sqrt{1 - (r')^2} \,\cosh \varphi(u) \, du,\\
\vspace{2mm}
x_4(u) =\int  \sqrt{1 - (r')^2} \,\sinh \varphi(u) \, du.
\end{array}
\end{equation}
Then the spacelike curve $c: \widetilde{z}(u) = \left( r(u), x_2(u), 0, x_4(u) \right)$ is a generating curve of a
constant mean curvature rotational surface of hyperbolic type.

Conversely, any constant mean curvature rotational surface of hyperbolic type is locally described by one of the cases given above.
\end{thm}

The proof of the theorem is similar to the proof of Theorem \ref{T:cmc elliptic}.

\vskip 2mm
In the case $r r'' + (r')^2 - 1 = 0$, i.e.
$r(u) = \pm \sqrt{u^2+2au+b},\; a=const\neq 0,\; b=const$,
the function  $\varphi(u)$ is expressed by the formula
$$\varphi(u)= \ds{\frac{2 \eta C}{\sqrt{\varepsilon(a^2-b)}}} \left(\frac{u+a}{2}\sqrt{u^2+2au+b} - \frac{\varepsilon(a^2-b)}{2} \ln \left|u+a+\sqrt{u^2+2au+b}\right| + d\right).$$

\vskip 3mm
\subsection{Constant mean curvature  rotational surfaces of parabolic type}

Now we shall  consider the rotational surface of parabolic type $\mathcal{M}'''$,  defined by formula \eqref{E:Eq-3}. 
The length of the mean curvature vector field of $\mathcal{M}'''$ is given by
\begin{equation*}  \label{E:Eq-8b}
\langle H, H \rangle = \ds{\frac{1}{4 f^2 f'^2} \left( f^2(x_1'' f' - x_1' f'')^2 - (f f'' + (f')^2)^2 \right)}.
\end{equation*}

In the case  $x_1'' f' - x_1' f'' = 0,\; f f'' + (f')^2 \neq 0$ the surface  $\mathcal{M}'''$ lies in a hyperplane  of $\E^4_2$. So,  we assume that $x_1'' f' - x_1' f'' \neq 0$ in an open interval $I \subset J$. 

In the following theorem we give a local description of constant mean curvature rotational surfaces of parabolic  type.

\begin{thm}\label{T:cmc parabolic}
Given a smooth function $f(u): I \subset \R \rightarrow \R$,
define the functions
$$\varphi(u) = \ds{ f' \left( A \pm \int \frac{1}{f'} \sqrt{\left((\ln|ff'|)'\right)^2 \pm 4 C^2}\; du \right)}, \quad \,\, C = const \neq 0, \; A = const,$$
and
$$x_1(u) = \ds{\int\varphi(u) du}; \qquad g(u) = \ds{\int\frac{\varphi^2(u) - 1}{2 f'(u)} \,du}.$$
Then the curve $c: \widetilde{z}(u) =  x_1(u)\, e_1 + f(u) \, \xi_1 + g(u) \, \xi_2$ is a spacelike curve generating a
constant mean curvature rotational surface of parabolic  type.

Conversely, any constant mean curvature rotational surface of parabolic type is locally constructed as described above.
\end{thm}

In the special case  $ff'' + (f')^2 = 0$, i.e. $f(u) = \pm \sqrt{2au+b}, \; a=const\neq 0, b=const$,  we obtain the following expression for the function $\varphi(u)$:
$$\varphi(u)= \ds{\frac{1}{\sqrt{2au+b}}} \left(A \pm\frac{2C B}{3a} \left(\sqrt{2au+b}\right) ^3\right), \; A = const, B = const \neq 0.$$

\vskip 5mm \textbf{Acknowledgments:}
The second author is partially supported by the National Science Fund,
Ministry of Education and Science of Bulgaria under contract
DFNI-I 02/14.

\end{document}